\documentclass[12pt]{amsart}
\usepackage{amssymb,amsthm} 
\usepackage{pxfonts}   

\newtheorem{thm}{Theorem}

\newtheorem{prop}{Proposition}
\newtheorem{lemma}{Lemma}
\theoremstyle{definition}

\newtheorem*{claim}{Claim}
\newtheorem*{example}{Example}
\theoremstyle{remark}
\newtheorem*{ack}{Acknowledgement}

\def\R{\mathbb{R}}
\def\C{\mathbb{C}}
\def\Z{\mathbb{Z}}
\def\N{\mathbb{N}}

\renewcommand{\P}{\mathbb{P}}
\def\eps{\varepsilon}

\def\cU{\mathcal{U}}

\DeclareMathOperator{\interior}{int}

\renewcommand{\epsilon}{\varepsilon}
\renewcommand{\setminus}{\smallsetminus}
\renewcommand{\emptyset}{\varnothing}
\def\SL{\mathrm{SL}}
\def\GL{\mathrm{GL}}

\begin{document}

\title[Cocycles over a minimal base]
{A uniform dichotomy for generic $\SL(2,\R)$ cocycles over a minimal base} 




\author[A.~Avila]{Artur Avila}
\address{CNRS UMR 7599,
Laboratoire de Probabilit\'es et Mod\`eles al\'eatoires.
Universit\'e Pierre et Marie Curie--Bo\^\i te courrier 188.
75252--Paris Cedex 05, France}
\urladdr{www.proba.jussieu.fr/pageperso/artur/}
\email{artur@ccr.jussieu.fr}

\author[J.~Bochi]{Jairo Bochi}
\address{Instituto de Matem\'atica -- UFRGS -- Porto Alegre, Brazil}
\urladdr{www.mat.ufrgs.br/$\sim$jairo}
\email{jairo@mat.ufrgs.br}

\begin{abstract}
We consider continuous $\SL(2,\R)$-cocycles over a minimal
homeomorphism of a compact set $K$ of finite dimension.
We show that the generic cocycle either is uniformly hyperbolic
or has uniform subexponential growth.
\end{abstract}

\date{November, 2006}

\maketitle

\section{Introduction}

In this paper we will consider $\SL(2,\R)$-valued cocycles over a minimal
homeomorphism $f:K \to K$ of a compact set $K$.
Such a cocycle can be defined as a pair $(f,A)$ where $A:K \to \SL(2,\R)$
is continuous.  The cocycle acts on $K \times \R^2$ by $(x,y) \mapsto
(f(x),A(x) \cdot y)$.  The iterates of the cocycle are denoted
$(f,A)^n=(f^n,A_n)$.

We say that $(f,A)$ is \emph{uniformly hyperbolic}
if there exists $\epsilon>0$ and $N>0$ such that
$\|A_n(x)\| \geq e^{\epsilon n}$ for every $x \in K$, $n \geq N$.
This is equivalent
to the existence of a continuous invariant splitting
$\R^2=E^u(x) \oplus E^s(x)$ such that vectors in $E^s(x)$ are
exponentially contracted by forward iteration and vectors in
$E^u(x)$ are exponentially contracted by backwards iteration --
see~\cite[proposition 2]{yoccoz}.

We say that $(f,A)$ has \emph{uniform subexponential growth} if for every
$\epsilon>0$ there exists $N>0$ such that
$\|A_n(x)\| \leq e^{\epsilon n}$ for every $x \in K$, $n \geq N$.
This condition is equivalent to the vanishing of the Lyapunov exponent for
all $f$-invariant probability measures (see propostion~\ref{p.subexp} below).
We recall that the Lyapunov exponent of the cocycle $(f,A)$
with respect to an $f$-invariant probability measure $\mu$
is defined as $L(f,A,\mu)=\lim \frac {1} {n} \int_K \log \|A_n\| \; d\mu$.

We say that a compact set $K$ has \emph{finite dimension} if it
is homeomorphic to a subset of some $\R^n$.
For instance, compact subsets of
manifolds (assumed as usual to be Hausdorff and
second countable) have finite dimension.
(For definitions of dimension and results concerning embedding in $\R^n$,
see e.g.~\cite{dimension}.)

\begin{thm} \label{generic uniformity}
Let $f:K \to K$ be a minimal homeomorphism of a compact set of finite
dimension.  Then for generic continuous
$A:K \to \SL(2,\R)$, either $(f,A)$ is uniformly hyperbolic or
$(f,A)$ has uniform subexpoential growth.
\end{thm}

In the case where $f$ is a minimal uniquely ergodic homeomorphism,
theorem \ref {generic uniformity} is contained in \cite {bochi}, which
shows that if $f:K \to K$ is a homeomorphism, then
for any ergodic $f$-invariant probability $\mu$,
the generic continuous $A:K \to \SL(2,\R)$ is such that
either $L(f,A,\mu)=0$ or
the restriction of $(f,A)$ to the support of $\mu$ is uniformly hyperbolic.
In general a minimal
homeomorphism may admit uncountably many ergodic invariant probability
measures, and in this paper we show how to treat them all at the same
time.  In the perturbation arguments, one often allows for loss of control
in some sets, and this must be compensated by showing that such sets can be
selected small.
The notion of smallness needed for our purposes involves simultaneous conditions on
all $f$-invariant probability measures.

\smallskip

It is unclear whether the hypothesis that $K$ has finite dimension is
actually necessary.  The following example shows that the minimality
hypothesis cannot be significantly weakened.

\begin{example}
Let $S^3$ be identified with $\C^2_*/\R_+$ (here $\C^2_*=\C^2 \setminus \{(0,0)\}$
and $\R_+=\{r \in \R;\; r>0\}$).
Fix some $\alpha \in \R$ with $\alpha/2\pi$ irrational and
let $f:S^3 \to S^3$ be given by
$$
(z,w) \mapsto (e^{i \alpha} z,e^{i \alpha}(z+w)) \bmod{\R_+} \, .
$$
Notice that if $h:S^3 \to \C\P^1$ with $h(z,w)=w/z$ is the usual Hopf fibration,
and $g:\C\P^1\to \C\P^1$ is $g(w)=w+1$, then $h\circ f = g \circ h$.
So $f$ has a unique minimal set, namely
$S=h^{-1}(\infty)$, where $f$ acts as a irrational rotation.
Let $A:S^3 \to \SL(2,\R)$ be any continuous map whose restriction to $S$
is given by
$$
(0,r e^{i \theta}) \mapsto
\begin {pmatrix} \cos (\theta+\alpha)&-\sin (\theta+\alpha)
\\ \sin (\theta+\alpha)&\cos (\theta+\alpha)\end{pmatrix}
\begin {pmatrix} 2&0\\0&1/2 \end{pmatrix}
\begin {pmatrix} \cos \theta&\sin \theta
\\ -\sin \theta&\cos \theta\end{pmatrix}.
$$
Then $(f|S,A|S)$ is uniformly hyperbolic: the associated splitting
is such that $E^u$ and $E^s$ are orthogonal and $E^u(0,r e^{i \theta})$
is generated by $\begin{pmatrix} \cos \theta \\
\sin \theta \end{pmatrix}$.
This splitting is topologically nontrivial, and hence cannot
extend to the whole $S^3$.  It follows that $(f,A)$ is not
uniformly hyperbolic and does not have uniform subexponential growth,
and the same properties hold for any small perturbation of $A$.
\end{example}

Using ideas from~\cite{bochi viana, cong}, one can adapt
the arguments of this paper to deal with $\GL(d,\R)$-valued cocycles.
The conclusion is that \emph{for a generic continuous $A:K \to \GL(d,\R)$
and for every $f$-invariant probability measure $\mu$,
the Oseledets splitting relative to $\mu$ coincides almost everywhere
with the finest dominated splitting of $(f,A)$.}

\begin{ack}
This research was partially conducted during the period A.A.\ served as a Clay Research Fellow.
The paper was written while A.A.\ was visiting UFRGS partially supported by Procad/CAPES.
\end{ack}

\section{Uniform subexponential growth}

Here we prove an equivalence stated at the introduction:

\begin{prop}\label{p.subexp}
Let $f:K \to K$ be homeomorphism of a compact set $K$.
and $A: K \to \SL(2,\R)$ be a continuous map.
Then the following are equivalent:
\begin{enumerate}
\item[(a)] $(f,A)$ has uniform subexponential growth:
for every $\epsilon>0$ there exists $N>0$ such that
$\|A_n(x)\| \leq e^{\epsilon n}$ for every $x \in K$, $n \geq N$.

\item[(b)] for every $\epsilon>0$ there exists $n>0$ such that
$\|A_n(x)\| \leq e^{\epsilon n}$ for every $x \in K$.

\item[(c)] $L(f,A,\mu)=0$ for every $f$-invariant probability $\mu$.
\end{enumerate}
\end{prop}

In the case that $f$ is uniquely ergodic, 
the proposition follows from~\cite[Theorem~1]{furman}.

\begin{proof}[Proof of the proposition]
(a) $\Rightarrow$ (b) is trivial;
(b) $\Rightarrow$ (c) follows from the fact that
$L(f,A,\mu) = \inf_n \frac{1}{n} \int_K \log\|A_n\| \; d\mu$.

\smallskip

We are left to prove (c) $\Rightarrow$ (a).
Assume that (a) does not hold.  Then there exists a sequence $x_k \in K$,
$n_k \to \infty$ such that $\|A_{n_k}(x_k)\| \geq e^{\epsilon n_k}$.  Let
$\mu_k=\frac{1}{n_k} \sum_{j=0}^{n_k-1} \delta_{f^j(x_k)}$.  Passing through
a subsequence, we may assume that $\mu_k$ converges to $\mu$, which is
$f$-invariant.  We claim that $L(f,A,\mu) \geq
\epsilon$.

Let $\delta>0$ and $s\in\N$ be fixed.  
It is enough to show that $\int \log \|A_s\| d\mu \geq
(\epsilon-\delta) s$.
Let $m_k=\lfloor n_k/s \rfloor$.  
Let $\nu_k=\frac{1}{s m_k} \sum_{j=0}^{s m_k-1} \delta_{f^j(x_k)}$.  
Notice that $\nu_k \to \mu$.
It is clear that if $k$
is large then $\|A_{s m_k}(f^i(x_k))\| \geq e^{(\epsilon-\delta) s m_k}$ for
$0 \leq i \leq s-1$.  Then
\begin{align*}
\int \log \|A_s\| d\nu_k&=\frac {1} {s m_k}
\sum_{i=0}^{s-1} \sum_{j=0}^{m_k-1}
\log \|A_s(f^{j s+i}(x_k))\|\\
&\geq \frac {1} {s m_k} \sum_{i=0}^{s-1} \log \|A_{s m_k}(f^i(x_k))\|
\geq s (\epsilon-\delta).
\end{align*}
The result follows.
\end{proof}

\section{Perturbation along segments of orbits}\label{s.perturbation}

In this section we assume that
$f:K \to K$ is minimal with no periodic orbits and
$A: K \to \SL(2,\R)$ is a continuous map
such that $(f,A)$ is not uniformly hyperbolic, but there exists
an $f$-invariant measure such that $L(f,A,\mu)>0$.
The aim here is to establish lemma~\ref{l.perturb} (see below).

We begin with an adaptation of lemma~3.4 from~\cite{cong}:

\begin{lemma}\label{l.cong}
For every $\epsilon>0$, there exists a non-empty open set
$W \subset K$, and $m \in \N$
such that:
\begin{itemize}
\item $W$, $f(W)$, \ldots, $f^{m-1}(W)$ are disjoint;
\item for all $x \in W$ and any non-zero vectors $\mathbf{v}$, $\mathbf{w}$,
there exists $M_0$, \ldots, $M_{m-1} \in \SL(2,\R)$ such that
$\|M_j - A(f^j(x)) \| < \eps$ and $M_{m-1} \cdots M_0 (\mathbf{v})$
is collinear to $\mathbf{w}$.
\end{itemize}
\end{lemma}

\begin{proof}
By assumption, there exists a $f$-invariant measure $\mu$ with non-vanishing exponents.
We can assume $\mu$ is ergodic (and hence non-atomic).
Then \cite[lemma~3.4]{cong} gives $m \in \N$ and a set $W$ which has
all properties we want, except it is not necessarily open.
(Notice that \cite{cong} gives the perturbed matrices in $\GL(2,\R)$,
but that is trivial to remedy.)
Reducing $W$, we can assume it to be a point.
Using the continuity of $A$, it is easy to see that $W$ can then be
slightly enlarged to become open.
\end{proof}

\begin{lemma}\label{l.perturb}
Given $\eps>0$, there exists arbitrarily large $N \in \N$ such that
for every $x\in K$ there exist
$L_0$, \ldots, $L_{N-1} \in \SL(2,\R)$ satisfying
$\|L_j - A(f^j(x)) \| < \eps$ and $\|L_{N-1} \cdots L_0\| < e^{\eps N}$.
\end{lemma}

Lemma~\ref{l.perturb} is an adaptation of lemma~5.1 from \cite{bochi}.
We will prove it using lemma~\ref{l.cong}.


If $A \in \SL(2,\R) \setminus \mathrm{SO}(2,\R)$,
then let $\mathbf{u}_A$ and $\mathbf{s}_A$ be unit vectors such that
$$
\| A (\mathbf{u}_A) \| = \|A\| \quad \text{and} \quad
\| A (\mathbf{s}_A) \| = \|A\|^{-1} \, .
$$
($\|\mathord{\cdot}\|$ always denotes euclidian norm.)
These vectors are unique modulo multiplication by $-1$.
Notice that $\mathbf{u}_A$ is orthogonal to $\mathbf{s}_A$,
and $A(\mathbf{u}_A)$ is collinear to $\mathbf{s}_{A^{-1}}$.


\begin{proof}[Proof of lemma~\ref{l.perturb}]
Let $W$ and $m$ be given by lemma~\ref{l.cong}.
Since $W$ is open and $f$ is minimal, there exists
$m_1 \in \N$ such that
\begin{equation}\label{e.m1}
\bigcup_{j=0}^{m_1} f^j(W) = K.
\end{equation}
Let $C > \eps + \sup_x \|A(x)\|$.
Take any $N \in \N$ be such that
$$
C^{4m_1 + 1} < e^{\eps N}/\sqrt 2  \, .
$$

From now on, let $x \in K$ be fixed.
We will explain how to find matrices $L_j$'s as in the statement of the lemma.
Let
$$
\Delta_j = \frac{\|A_j(x)\|}{\|A_{N - j}(f^j(x))\|} \, , \quad
\text{for $0 \le j \le N$.}
$$
Then $\Delta_{N} = 1/\Delta_0$ and
$C^{-2} < \Delta_{j+1} / \Delta_j < C^2$.
It follows that there exists $j_0$ such that
$C^{-1} < \Delta_{j_0} < C$.

Due to~\eqref{e.m1}, there exists $j_1$ such that
$j_0 \le j_1 \le j_0 + m_1$ and $f^{j_1}(x) \in W$.
We can assume that $j_1 + m \le N$, because
otherwise we would have $N < j_0 + 2m_1$, so
$$
\|A_{N}(x)\| = \Delta_{N} \le \Delta_{j_0} C^{4m_1} \le C^{4m_1+1} < e^{\eps N},
$$
and then there would be nothing to prove.

Let
$X = A_{j_1}(x)$ 
and $Z = A_{N - j_1 -m}(f^{j_1+m}(x))$.
Let $M_0$, \ldots, $M_{m-1}$ be given by lemma~\ref{l.cong}
so that $\|M_i - A(f^{j_1 + i}(x))\|< \eps$
and $M_{m-1} \cdots M_0 (\mathbf{s}_{X^{-1}})$ is collinear to $\mathbf{s}_Z$.

For $0 \le j < N$, let
$$
L_j =
\begin{cases}
M_{j-j_1} &\text{if $j_1 \le j < j_1 + m$,}\\
A(f^j(x)) &\text{otherwise.}
\end{cases}
$$

Write $Y = M_{m-1} \cdots M_0$, so $L_{N-1} \cdots L_0 = ZYX$.
Notice $YX \cdot \mathbf{u}_X$ is collinear to $\mathbf{s}_Z$.
So
\begin{align*}
\|ZYX (\mathbf{u}_X)\|
&\le \|A_{j_0}(x) (\mathbf{u}_X)\| \cdot C^{j_1-j_0}\cdot\|Y\|\cdot\|Z(\mathbf{s}_Z)\| \\
&\le \|A_{j_0}(x)\| \cdot C^{m_1+m} \cdot \|Z\|^{-1} \\
&\le C^{2m_1+2m} \|A_{j_0}(x)\| \cdot \|A_{N-j_0}(f^{j_0}(x))\|^{-1} \\
&= C^{2m_1+2m} \Delta_{j_0}
< C^{4m_1+1}. \\
\intertext{Also,}
\|ZYX (\mathbf{s}_X)\|
&\le \|X\|^{-1} \cdot \|Y\| \cdot \|Z\| \\
&\le \|A_{j_0}(x)\|^{-1} C^{m_1} \cdot C^m \cdot C^{m_1} \|A_{N-j_0}(f^{j_0}(x))\| \\
&= C^{2m_1+m} \Delta_{j_0}^{-1}
< C^{3m_1+1}.
\end{align*}
We have shown that
$\max\big(\|ZYX(\mathbf{u}_X)\|, \|ZYX(\mathbf{s}_X)\|\big)<
C^{4m_1+1} < e^{\eps N}/\sqrt 2$.
Since $\mathbf{u}_X \perp \mathbf{s}_X$, we conclude that
$\|ZYX\| < e^{\eps N}$, as wanted.
\end{proof}

\section{Tiling $K$}

A Borel set ${X \subset K}$ is said to be a \emph{zero probability set}
for the homeomorphism $f:{K \to K}$ if $\mu(X)=0$ for every
$f$-invariant probability measure $\mu$.

Our goal in this section is to prove lemma~\ref{l.boundary} below.

\begin{lemma} \label{l.boundary}
Let $f:K \to K$ be a homeomorphism of a compact set of finite dimension
with no periodic orbits.
There exists a basis of the topology of $K$ consisting of sets $U$
such that $\partial U$ is a zero probability set.
\end{lemma}

To prove lemma~\ref{l.boundary}, we will need lemmas~\ref{l.extension}
and \ref{l.periodic}.

\begin{lemma}\label{l.extension}
Let $f:K \to K$ be a homeomorphism of a compact set of finite dimension.
Then there exists $d>0$, an embedding $s:K \to \R^d$ and a homeomorphism
$g:\R^d \to \R^d$ such that $s \circ f=g \circ s$.
\end{lemma}

\begin{proof}
This result is well known, but we reproduce the proof for convenience.
We may assume that $K \subset \R^n$ for some $n$.  Let $\phi:\R^n \to \R^n$
and $\psi:\R^n \to \R^n$ be continuous extensions of $f$ and $f^{-1}$, respectively.
Let $d=2n$, $s(x)=(x,f(x))$ and $g(x,y)=(y,{x+\phi(y)-\psi(y)})$.
\end{proof}

\begin{lemma} \label{l.periodic}
Let $f:\R^d \to \R^d$ be a homeomorphism.
and let $S$ be a compact manifold of dimension $d-1$.
For a generic continuous $\psi:S \to \R^d$, and for every sequence of integers
$j_0<\cdots<j_d$, $\bigcap_{k=0}^d f^{j_k}(\psi(S))$ is contained in the set
of periodic orbits of $f$.
\end{lemma}

\begin{proof}
Let $g_i: V_i \to \R^d$ be a countable family of charts
so that there are exists a basis of the topology of $S$ formed by sets
$U_i$ such that $\overline{U_i} \subset V_i$.
Given sequences $\mathbf{j} = (j_0<\cdots<j_d)$
and $\mathbf{i} = (i_0,\ldots,i_d)$ such that
$\overline{U_{i_0}}$, \ldots, $\overline{U_{i_d}}$ are disjoint,
let $\cU_{\mathbf{j},\mathbf{i}}$
be the set of continuous maps $\psi: S \to \R^d$ such that
$$
\bigcap\nolimits_{k=0}^d f^{j_k} \circ \psi \big(\overline{U_{i_k}}\big)
= \emptyset.
$$

\begin{claim}
The set $\cU_{\mathbf{j},\mathbf{i}}$ is open and dense in $C(S,\R^d)$.
\end{claim}

Assuming the claim for the moment, let us conclude the proof of the lemma.
Let $\psi$ belong to the residual set $\bigcap\cU_{\mathbf{j},\mathbf{i}}$.
Assume that $z \in \bigcap_{k=0}^d f^{j_k}(\psi(S))$
for some $j_0<\cdots<j_d$.
Let $x_k \in S$ be such that $f^{j_k}(\psi(x_k))=z$, for $0\le k \le d$.
If all $x_k$ where distinct then we could take chart domains
$U_{i_k} \ni x_k$ with
$\overline{U_{i_0}}$, \ldots, $\overline{U_{i_d}}$ disjoint.
This would contradict $\psi \in \cU_{\mathbf{j},\mathbf{i}}$.
We conclude that at least two $x_k$'s coincide.
It follows that $z$ is a periodic point of $f$.

\smallskip

Now, the proof of the claim.
Openness is clear; it remains to show density.
For simplicity of writing, let $i_k=k$.
Reducing the sets $V_k$ if necessary, we can assume they are disjoint
(but still with $V_k \supset \overline{U_k}$).
By a small perturbation of $\psi$ supported on $V_k$, we may assume that
$f^{j_k} \circ \psi \circ g_k^{-1}$ is smooth in a neighborhood of $g_k(U_k)$.
In other words, letting
$L = \overline{g_0(U_0)} \times \cdots \times \overline{g_d(U_d)}
\subset (\R^{d-1})^{d+1}$,
there is a neighborhood $W \supset L$ such that the map
$$
G : W \to (\R^d)^{d+1}, \quad
G=(f^{j_0} \circ \psi \circ g_0^{-1},\ldots,f^{j_d} \circ \psi \circ g_d^{-1})
$$
is smooth.
Perturbing $\psi$ again, we may assume that
$G$ is transverse to the diagonal
$D=\big\{(x,\ldots,x) \in (\R^d)^{d+1};\; x \in \R^d\big\}$ at $L$.
Since the diagonal has codimension $d^2$ and $W$ has dimension $d^2-1$, this
implies that $G(L)$ does not intersect $D$.
That is, $\bigcap_{k=0}^d f^{j_k}\circ \psi \big(\overline{U_k} \big) =\emptyset$.
\end{proof}

\begin{proof}[Proof of lemma \ref {l.boundary}]
By lemma~\ref{l.extension}, we can assume that $K \subset \R^d$ and
$f:K \to K$ is the restriction of a homeomorphism $f:\R^d \to \R^d$.

Let $x_0 \in K$, and let $\epsilon>0$.
We need to show that there exists an
open set $U \subset \R^d$ containing $x_0$ and of diameter at most
$4 \epsilon$, such that
$\mu(\partial U)=0$ for every $f$-invariant probability $\mu$ supported on $K$.

Let $B$ be the closed unit ball in $\R^d$, and let $S=\partial B$.
Let $\phi:B \to \R^d$ be given by $\phi(y)=x_0+\epsilon y$.
By lemma~\ref{l.periodic},
there exists a continuous map $\psi:S \to \R^d$ such that
$\|\psi(y)-\phi(y)\|<\epsilon$ for every $y \in S$ and
such that for every $x \in K$, the set of $\ell \in \Z$ such that
$f^{-\ell}(x) \in \psi(S)$ has cardinality at most $d+1$.  By the Ergodic
Theorem, $\mu(\psi(S))=0$ for every $f$-invariant probability measure.
Let $U$ be the
connected component of $x_0$ in $\R^d \setminus \psi(S)$.
Then the diameter of $U$ is at most $4 \epsilon$
and $\partial U \subset \psi(S)$, so $U$ has all the desired
properties.
\end{proof}

\section{Proof of theorem \ref {generic uniformity}}

Let $f:K \to K$ be a minimal homeomorphism of a compact set of finite dimension.
From now on, we assume that $f$ does not have periodic points;
otherwise the result is obvious.

The idea of the following lemma and its proof come from \cite{rokhlin}.

\begin{lemma}\label{l.castle}
For any $N\in \N$, there exists an open set $B \subset K$ such that:
\begin{itemize}
\item the return time from $B$ to itself
via $f$ assumes the values $N$ and $N+1$ only;
\item $\partial B$ has zero probability.
\end{itemize}
\end{lemma}

\begin{proof}
Given $N$, 
there exists $n_1 \in \N$ such that
if $n \ge n_1$ then there are $\ell$, $\ell' \ge 0$ satisfying
$n = \ell N + \ell' (N+1)$.
By lemma~\ref{l.boundary} and the fact that $f$ has no periodic
points, we can take an open set $U$
such that $\partial U$ has zero probability and
$U$, $f(U)$, \ldots, $f^{n_1}(U)$ are disjoint.
Then, for some $n_1> n_1$, $\bigcup_{n= 0}^{n_1} f^n (U) = K$.
For each $n \in N$, let
$$
U_n = U \cap f^{-n}(U) \setminus \bigcup_{j=1}^{n-1} f^{-j}(U),
$$
that is, the set of points $x \in U$ such that the first return time to $U$ is $n$.
Notice that $\partial U_n$ has zero probability.
For each $n$ between $n_1$ and $n_1$ for which $U_n \neq \emptyset$,
consider the tower of height $n$ with base $U_n$.
Then break this tower into towers of heights $N$ or $N+1$.
In this way we cover $K$ with finitely many towers of heights $N$ or $N+1$.
Let $B$ be the union of the bases of the towers.
Since $B$ is a finite union of iterates of $U_n$'s, we have that
$\partial B$ has zero probability.
\end{proof}

\begin{lemma}\label{l.freq}
If $L$ is a compact set with zero probability then
for every $\eps>0$, there exists an open set $V \supset L$
and $n_0 \in \N$ such that
\begin{equation}\label{e.freq}
\frac{1}{n} \; \#\{j; \; 0 \le j \le n-1, \ f^j(x) \in V\} < \eps
\quad \text{for all $x\in K$, $n \ge n_0$.}
\end{equation}
\end{lemma}

\begin{proof}
Otherwise there exist $\epsilon>0$,
$x_k \in K$, open sets $V_k \subset K$ and integers $n_k \to \infty$
such that $\overline V_{k+1} \subset \interior V_k$,
$\bigcap V_k=L$ and $\mu_k(V_k)>\epsilon$, where
$\mu_k=\frac{1}{n_k} \sum_{j=0}^{n_k-1} \delta_{f^j(x_k)}$ (here
$\delta_x$ is the Dirac mass on $x$).
If $\mu$ is a weak-$*$ limit of the probability measures
$\mu_k$ then $\mu$ is $f$-invariant.  Moreover, we have
$\mu(V_k) \geq \liminf \mu_j(V_j) \geq \epsilon$ for every $k$.
We conclude that $\mu(L) \geq \epsilon$, contradiction.
\end{proof}

\begin{proof}[Proof of theorem~\ref{generic uniformity}]
The set $\mathit{UH}$
of continuous $A:K \to \SL(2,\R)$ such that $(f,A)$ is uniformly hyperbolic
is open.
For any $\eps>0$, let $\cU_\eps$ be the set of $A$'s such that
there is $n_1 \in \N$ such that
$\frac{1}{n_1} \log \| A_{n_1}(x) \| < \eps$
for all $x \in K$;
then $\cU_\eps$ is also an open set.
If $A \in \bigcap_{\eps>0} \cU_\eps$ then $(f,A)$ has
uniform subexponential growth, by proposition~\ref{p.subexp}.
So to prove the theorem it suffices to show that $\cU_\eps$ is dense
in the complement of $\mathit{UH}$.

\smallskip

Fix $\eps>0$ and $A \notin \mathit{UH}$.
We can assume that $(f,A)$ has positive exponent for some
invariant probability, otherwise (by proposition~\ref{p.subexp})
$f \in \cU_\eps$ already.
We will find $\tilde A$ close to $A$ such that $\tilde A \in \cU_\eps$.

Let $c>0$ be such that $\|A\|+ \eps < e^c$.
Let $N \in \N$ be given by lemma~\ref{l.perturb} (notice hypotheses from \S\ref{s.perturbation} hold).
We can assume $\eps N> c$.
Then let $B \subset K$ be the set given by lemma~\ref{l.castle}.
Let $\delta>0$ be small enough so that
$$
0\le j \le N, \ d(x,y)<\delta \ \Rightarrow \ \| A(f^j(x)) - A(f^j(y)) \| < \eps.
$$

Cover the closure of $B$ by open sets $W_1$, \ldots, $W_k$
with diameter less than $\delta$.
By lemma~\ref{l.boundary}, we can take each $W_i$ so that
$\partial W_i$ has zero probability.
Let $U_i = W_i \setminus \bigcup_{j<i} \overline {W_j}$.

Let $B_\ell$ be the set of points in $B$ whose first return to $B$
occurs in time $\ell$; then $B = B_N \sqcup B_{N+1}$.
Let
$$
L = \bigcup_{\ell=N}^{N+1} \bigcup_{i=0}^{\ell-1}
\partial (B_\ell \cap U_i) \, .
$$
Then $L$ has zero probability.
By lemma~\ref{l.freq}, there exists an open set $V \supset L$
and $n_0 \in \N$ such that
\begin{equation}\label{e.freq2}
\frac{1}{n} \; \#\{j; \; 0 \le j \le n-1, \ f^j(x) \in V\} < \frac{\eps}{N+1}
\quad \text{for all $x\in K$, $n \ge n_0$.}
\end{equation}

For each $\ell= N$ or $N+1$ and $i=1$, \ldots, $k$,
we choose a point $x_{\ell,i}$ in $B_\ell \cap U_i$.
Applying lemma~\ref{l.perturb},
we find $L_{\ell,i,0}$, \ldots, $L_{\ell,i,N-1}$
so that
$$
\|L_{\ell,i,j} -  A(f^j(x_{\ell,i})) \| < \eps \ \forall j=0,\ldots, \ell -1
$$
and $\| L_{\ell,i,N-1} \cdots L_{\ell,i,0} \| < e^{\eps N}$.
In the case $\ell =N+1$, define also $L_{\ell,i,N} = A(f^j(x_{\ell,i}))$.
So for any $\ell$, $\| L_{\ell,i,\ell-1} \cdots L_{\ell,i,0} \| < e^{2 \eps \ell}$.

We want to define a map $\tilde A: K \to \SL(2,\R)$.
Begin defining
$$
\tilde A = L_{\ell,i,j}
\text{ on the set $f^j\big(B_\ell \cap U_i \setminus V\big)$,}
$$
where $N\le\ell\le N+1$, $1\le i \le k$, and $0 \le j \le \ell-1$.
(Notice these sets are disjoint.)
It remains to define $\tilde A$ on the rest of $K$.
The map $A^{-1} \tilde{A}$ is already defined
on a compact subset of $K$, and takes values on the
$(e^c+1)\eps$-neighborhood of $\mathrm{Id}$ in $\SL(2,\R)$.
That domain is homeomorphic to $\R^3$ (provided $\eps$ is not too big),
so by the Tietze's extension theorem we can continuously extend
$A^{-1} \tilde{A}$ to the whole $K$.
So we obtain $\tilde A: K \to \SL(2,\R)$ with
$\|\tilde{A} -A \| < e^c(e^c+1)\eps$.

\smallskip

Now let $n > \max(n_0, (N+1) / \eps)$.
Fix $x\in K$.
We will give an upper bound for $\|\tilde{A}_n(x)\|$.
Write
$$
n = p + \ell_1 + \ell_2 + \cdots + \ell_r + q
$$
in a way such that
$0 \le p, q \le N+1$, $N \le \ell_i \le N+1$,
and the points
$$
x_1 = f^p(x), x_2 = f^{p+\ell_1}(x), \ldots, x_r = f^{p+\ell_1+ \cdots + \ell_r}(x)
$$
are exactly the points in the segment of orbit
$x$, $f(x)$, \ldots, $f^{n-1}(x)$ that belong to $B$.
We have
$$
\|A_n(x)\| \le e^{2c(N+1)} \prod_{i=1}^r \|A_{\ell_i}(x_i)\| \, .
$$
We will use the bounds:
$$
\|A_{\ell_i}(x_i)\| \le
\begin{cases}
e^{2\eps \ell_i} &\text{if $x_i \in B \setminus V$,}\\
e^{c (N+1)}     &\text{if $x_i \in V$.}
\end{cases}
$$
By~\eqref{e.freq2},
there are at most $(\eps/(N+1)) n$ points $x_i$ that belong to $V$.
So
$$
\|A_n(x)\| \le e^{2c(N+1)} \cdot (e^{c(N+1)})^{(\eps/(N+1)) n} \cdot e^{2\eps \sum \ell_i}
           < e^{(3c+2)\eps n}.
$$
This proves that $\tilde{A} \in \cU_{(3c+2)\eps}$.
\end{proof}



\begin{thebibliography}{BV}

\bibitem[B]{bochi}
J.~Bochi.
Genericity of zero Lyapunov exponents.
\textit{Ergodic Theory Dynam.\ Systems}  22  (2002), 1667--1696.

\bibitem[BV]{bochi viana}
J.~Bochi and M.~Viana.
The Lyapunov exponents of generic volume preserving and symplectic maps.
\textit{Annals of Math.\ }161 (2005), 1423--1485.

\bibitem[C]{cong}
Nguyen Dinh Cong.
A generic bounded linear cocycle has simple Lyapunov spectrum.
\textit{Ergodic Theory Dynam.\ Systems}  25  (2005), 1775--1797.

\bibitem[EP]{rokhlin}
S.~Eigen, V.~S.~Prasad.
Multiple Rokhlin Tower Theorem: A Simple Proof.
\textit{New York J.\ Math.\ }3A (1997/98), 11--14.

\bibitem[F]{furman}
A.~Furman.
On the multiplicative ergodic theorem for uniquely ergodic systems.
\textit{Ann.\ Inst.\ H.\ Poincar\'{e} Probab.\ Statist.\ }33 (1997), 797--815.

\bibitem[N]{dimension}
J.-I.~Nagata.
\textit{Modern Dimension Theory}.
Wiley, 1965.

\bibitem[Y]{yoccoz}
J.-C.~Yoccoz.
Some questions and remarks about $SL(2,\mathbb{R})$ cocycles.
\textit{Modern dynamical systems and applications}, 447--458,
Cambridge Univ.~Press, 2004.


\end{thebibliography}
\end{document}